\documentclass{amsart}
\usepackage{amsmath,amssymb}
\def\p1i{\pi_1^{\infty}}
\let\dfr=\rightarrow

\let\ =\quad

\let\ti=\widetilde
\newtheorem{theo} {Theorem}
\newtheorem{defi}{Definition} 
\newtheorem{cor}{Corollary} 
\newtheorem{lem}{Lemma} 
\newtheorem{prop}{Proposition} 
\newtheorem{rem}{Remark}

\title{A refinement of the simple connectivity at infinity of groups}

\author{Louis Funar} 
\address{Institut Fourier BP 74, UMR 5582,  Universit\'e Grenoble I, 38402 
Saint-Martin-d'H\`eres Cedex, France } 
\email{funar@fourier.ujf-grenoble.fr}  
\author{Daniele Ettore Otera} 
\address{Universit\'e Paris-Sud, Math\'ematiques, B\^at 425, 91405 Orsay Cedex, France} 
\email{Daniele.Otera@math.u-psud.fr} 
\date{February 26, 2003.}
 
\begin{document} 
 
%\maketitle 
 
\begin{abstract} 
We give another proof for a result of Brick (\cite{Br}) stating that 
the simple connectivity at infinity is a geometric
property of finitely presented groups. This allows us to define the 
rate of vanishing of $\p1i$ for those groups which are simply connected 
at infinity. Further we show that this rate is linear for cocompact
lattices in nilpotent and semi-simple Lie groups, and in particular for
fundamental groups of geometric 3-manifolds. 
 
\vspace{0.1cm}
 
\vspace{0.1cm}
\noindent {\bf MSC Subject:} 20 F 32, 57 M 50.  
\end{abstract}

\maketitle

\section{Introduction}

 The first aim of this note is to prove the quasi-isometry 
 invariance of the simple connectivity at infinity for groups, in contrast with  
 the case of spaces. We recall that:

\begin{defi} 
 The metric spaces $(X,d_X)$ and $(Y,d_Y)$ are quasi-isometric 
 if there are constants $\lambda$, $C$ and maps $f: X \dfr Y$, $g: Y \dfr X$ (called $(\lambda ,C)$-quasi-isometries) 
 such that  the following: 
$$d_Y (f(x_1 ),f(x_2 ))\leqslant \lambda d_X (x_1 ,x_2)+C,$$ 
$$d_X (g(y_1 ),g(y_2 ))\leqslant \lambda d_Y (y_1 ,y_2)+C,$$ 
$$d_X(fg(x),x)\leqslant C,$$ 
$$d_Y(gf(y),y)\leqslant C,$$ 
hold true for all $x,x_1 ,x_2 \in X, y,y_1 ,y_2  \in Y$. 
\end{defi}

\begin{defi} 
 A  connected, locally compact, topological space $X$ with $\pi_1 X=0$ 
 is simply connected at infinity (abbreviated s.c.i. and one writes 
also $\pi_1^{\infty}X=0$) if for each compact $k\subseteq X$ there exists a larger compact 
$k\subseteq K\subseteq X$ such that any closed loop in $X - K$ is null homotopic in 
$X-k$. 
\end{defi} 
 \begin{rem}
The simple connectivity at infinity is not a quasi-isometry invariant  
of spaces (\cite{Ot}). In fact $(S^1 \times {\bf R}) \underset{S^1 \times \mathbb{Z} }{\cup} D^2$  
\ and \ $(S^1 \times {\bf R}) \underset{S^1 \times \{ 0 \} }{\cup} D^2$ are simply connected,  
quasi-isometric spaces although the first is simply connected at infinity while the second is not. 
\end{rem} 
\noindent This notion extends to a group-theoretical framework as follows (see \cite{BrM}, p.216): 
 
\begin{defi} 
A group $G$ is simply connected at infinity if for some (equivalently any)  finite complex $X$ 
such that $\pi_1 X=G$ one has $\p1i \ti{X}=0$, where $\ti{X}$ denotes the universal covering of $X$. 
\end{defi} 
The independence on the particular complex $X$ is proved in \cite{Ta} and \cite{Ot}. 
Roughly speaking the simple connectivity at infinity depends only on the 2-skeleton and  
any finite 2-complex corresponds to a presentation of $G$. Since  
Tietze transformations act transitively on the set of group presentations  
it suffices to check the invariance under such moves. We refer to \cite{Ta} for details.  
 
All groups considered in the sequel will be finitely generated and  a 
system of generators  
determines a word metric on the group. Although this depends on the
chosen  generating set  
the different word metrics are quasi-isometric. Therefore properties  which are  
invariant under quasi-isometries are independent on the particular word metric and will  
be called {\em geometric properties}.
It is well-known that being finitely presented, 
word hyperbolic or virtually free are geometric properties, while 
being virtually solvable or virtually torsion-free are not geometric
(see \cite{DE}). 

 Our main result is : 
 
\begin{theo}[\cite{Br}]
The simple connectivity at infinity of groups is a geometric property. 
\end{theo} 
 
\noindent This was originally proved by Brick in \cite{Br}.  
We provide  a simpler and more conceptual  proof, by analyzing the colored  Rips complex. 
 
\begin{rem} 
It seems that the fundamental group at infinity, whenever it is well-defined (see \cite{GM} for an extensive discussion on this topic), should also be a quasi-isometry invariant of the group. 
\end{rem} 

\begin{defi}
Let $X$  be a simply connected non-compact metric space with
$\p1i X=0$. The rate of vanishing of $\p1i$, denoted $V_X(r)$, 
is the infimal $N(r)$ with the property that any loop which sits 
outside the ball $B(N(r))$ of radius $N(r)$ 
bounds a 2-disk outside $B(r)$.
\end{defi}
\begin{rem}
It is easy to see that $V_X$ can be an arbitrary large function. 
\end{rem}
It is customary to introduce 
the following equivalence relation on functions: 
$f\sim g$  if 
there exists constants $c_i, C_j$ (with $c_1,c_2 >0$) such that 
\[ c_1 f(c_2 R)+ c_3 \leq g(R) \leq C_1 f(C_2 R) + C_3. \]
It is an easy consequence of the proof of theorem 1 that 
the equivalence class of $V_X(r)$ is a quasi-isometry invariant. 
In particular $V_G=V_{\widetilde{X_G}}$ 
is a quasi-isometry invariant of the group $G$, where
${\widetilde{X_G}}$ is  the universal covering space 
of a compact simplicial complex $X_G$, 
with $\pi_1(X_G)=G$ and $\p1i(G)=0$. 

\begin{rem}
For most groups $G$ coming from geometry $V_G$ is trivial, 
i.e. linear. Obviously  if $M$ has an Euclidean structure 
then $V_{\pi_1(M)}$ is linear. Since metric balls in the hyperbolic 
space are diffeomorphic to standard  balls in ${\bf R}^n$ 
one derives that  $V_{\pi_1(M)}$ is linear for any compact  
hyperbolic manifold $M$.  
\end{rem}
\begin{rem}
Notice that there exists (see \cite{D}) word hyperbolic groups $G$ 
(necessary of dimension $n\geq 4$ by \cite{BM}) 
which are not simply connected at infinity  and hence $V_G$ is not
defined. Moreover if $G$ is a word hyperbolic torsion-free group 
with $\p1i(G)=0$ then it seems that $V_G$ is linear. 
\end{rem}

\begin{theo}
$V_G$ is linear for uniform lattices in: 
\begin{enumerate}
\item semi-simple Lie groups. 
\item nilpotent groups. 
\item solvable stabilizers of horospheres in product of symmetric
  spaces of rank at least two, or generic horospheres in products of
  rank one symmetric spaces. 
\end{enumerate}
\end{theo}
Interesting examples of groups for which $\p1i(G)=0$ are 
the (infinite) fundamental groups of geometric 3-manifolds (and conjecturally 
of all 3-manifolds). We can show that: 

\begin{cor}
The fundamental 
groups of geometric 3-manifolds have linear $V_G$. 
\end{cor} 

\begin{rem}
The existence of  groups $G$ acting freely and cocompactly 
on ${\bf R}^n$, which have super-linear $V_G$ seems most likely.
The examples described in (\cite{Gr2}, section 4),
which have large acyclicity radius,  strongly support this claim. 
The first point is that the rate of vanishing
of $\p1i$ is rather related to higher (i.e. dimension $n-2$) 
connectivity radii, which are less understood. The second 
difficulty is that these groups are
not s.c.i. The simplest way to overcome 
it is to consider group extensions. For instance $\p1i(G\times {\bf Z}^2)=0$, 
for any finitely presented group $G$; alternatively 
$\p1i(V^n\times {\bf R})=0$ 
for any contractible manifold $V^n$ ($n\geq 2$).   
However this idea does not work because  
$V_{G\times {\bf Z}^2}$ is always linear. 
\end{rem} 
\begin{rem}
One needs some extra arguments in order to extend the proof to all 
solvable Lie groups. 
However, it seems very likely that 
cocompact lattices in  all connected Lie groups have linear 
rate of vanishing of $\p1i$.  
\end{rem}

\vspace{0.2cm}
 
\noindent $\bf{Acknowledgements}:$ The authors are indebted to
C. Dru\c tu, P. Pansu and F. Paulin for helpful conversations 
and to the referee for the useful comments improving the exposition.

\section{Proof of Theorem 1} 
 
\noindent For positive $d$, set $P_d(G)$ for the simplicial 
 complex defined as follows: 
\begin{itemize} 
\item its vertices are the elements of $G$,  
\item the elements $x_1......x_n$ of $G$ span an $n$-simplex, if  
$d(x_i,x_j)\leq d$ for all $i,j$ (where $d(. ,.)$ is the word metric). 
\end{itemize} 
 
\begin{rem} 
If $G$ is $\delta$-hyperbolic then $P_d(G)$ is contractible as soon as $d > 4\delta+1$ (see \cite{GdlH}). 
\end{rem}
Although $P_d(G)$ is not contractible in general, one can prove that 
it is simply connected under a mild restriction. 
Let $G=<x_1,...,x_n| R_1,...,R_p>$ be a presentation of $G$ and $r$ denotes the maximum length  
among the relators $R_i$.   
 
\begin{lem} 
If $2d > r$, then $\pi_1(P_d(G))=0$. 
\end{lem} 
\begin{proof} 
Let $l=[1,\gamma_1,\gamma_2,...,\gamma_n,1]$ be a  (simplicial) 
loop in $P_d(G)$ based at the identity. 
Two successive vertices of $l$ are at distance  at most $d$. 
One can interpolate between two consecutive $\gamma_j$'s   
a sequence of elements of $G$ (of length at most $d$), 
consecutive ones being adjacent when viewed 
as elements of the Cayley graph (hence at distance one). 
The product of elements corresponding to these edges of 
length one of $l$ is trivial in $G$. Therefore it is a product of 
conjugates of relators: $\prod g_i^{-1}R_ig_i$. 
The diameter of each $R_i$ is at most $r/2$, and the assumption  $2d>r$, implies that  
each loop $g_i^{-1}R_ig_i$ is contractible in $P_d(G)$. 
This ends the proof. 
\end{proof}

The natural group action of $G$ on itself by left translations gives
rise to an action on $P_d(G)$. In particular, if $G$ has no torsion
then it acts freely on $P_d(G)$  and $P_d(G)/G=X$ is a compact
simplicial complex with $\pi_1(P_d(G)/G)=G$.  
 
\begin{prop} 
The vanishing of $\p1i$ is a geometric property of torsion-free
groups. 
\end{prop} 
\begin{proof} 
One has to show that if the group  $H$ is quasi-isometric to  $G$ then
$\p1i P_d (G)=0$ implies that  $\p1i P_a (H)=0$ for large enough $a$. 
 
Let $f:H\to G$ and $g:G\to H$ be  $(k,C)$-quasi-isometries 
between $G$ and $H$. Fix $x_0 \in P_a(H)$ and 
$f (x_0)\in P_d(G)$ as base points. 

\begin{lem}
If $\p1i P_d (G)=0$ then $\p1i P_D (G)=0$ for $D\geq d$. 
\end{lem}
\begin{proof}
An edge in $P_D(G)$ corresponds to a path (of length uniformly bounded
by $\frac{D}{d}+1$) in $P_d (G)$. 
Thus a loop  $l$ in $P_D (G)$ at distance
at least $R$ from a given point corresponds to a loop $l'$ 
in $P_d (G)$ at distance at least $R-\frac{D}{d}-1$ from the same point. 
By assumption  $l'$ will bound a 2-disk $D^2$ far away in
$P_d(G)$. 

Now the union of an edge $[x_1x_n]$  in $P_D (G)$ 
and its corresponding path $[x_1,x_2,...,x_n]$ in $P_d(G)\subset P_D(G)$ form 
the boundary of a 2-disk in $P_D(G)$, which 
is triangulated by using the triangles $[x_1,x_j,x_n]$. 
Consider one such triangulated 2-disk for each edge of $l$  and 
glue  to the previously obtained  $D^2$ to get  a 2-disk 
in $P_D(G)$ bounding $l$ and far away.  
\end{proof}

By hypothesis for each $r$  there exists $N(r) >0 $ such that  every 
loop $l$ in $P_d(G)$ satisfying  
$d(l,f(x_0))> N(r)$ bounds a disk outside   
$B(f(x_0),r)$. This  means that there 
exists a simplicial map $\varphi:D^2\to P_d(G)-B(f(x_0),r)$ 
such that $\varphi(\partial D^2)$ is the
given loop $l$, when $D^2$ is suitably triangulated.  
A loop $l=[x_1,x_2,...,x_n,x_1]$ in $P_d(G)$, based at $x_1$, is the 
one-dimensional simplicial sub-complex with vertices $x_j$ and 
edges $[x_i x_{i+1}]$, $i=1,n$ (with the convention  $n+1=1$).

Set $M(R)=kN(kR+kC+3C)+3C$. We claim that:
\begin{lem}
Any loop $l$ in $P_a(H)$ sitting outside the ball $B(x_0,M(R))$ 
bounds a 2-disk not intersecting $B(x_0,R)$. 
\end{lem}
 \begin{proof}
Set $l=[x_1,...,x_n]$. Using the previous lemma one can assume that 
$d$ is large enough such that $\frac{d-C}{k} >1$. As in lemma 1 
one can add extra vertices between the consecutive ones such that  
$d(x_i,x_{i+1})\leq \varepsilon$ holds,  where $k\varepsilon + C = d$.

The image $f(l)=[f(x_1),....,f(x_n),f(x_1)]$ of the loop $l$ has the  
property that $d(f(x_i),f(x_{i+1}))\leq k\varepsilon +C$.   
 Using  $d(x,gf(x))\leq C$ one obtains that $d(gf(x),gf(y))\geq d(x,y)-2C$, which implies  
$d(x,y)\leq k d(f(x),f(y)) +3C$ and thus:  
\[d(f(x),f(y))\geq \frac{d(x,y)-3C}{k}, \mbox{ for all } x, y \in P_a(H).\] 
  
\noindent From this inequality one derives  that:   
\[d(f(x_i),f(x_0)\geq (M(R)-3C)/k=N(kR+kC+3C)\]  
and thus the  loop  $f(l)$ sits outside the ball 
$B(f(x_0), N(kR+kC+3C))$, and hence by assumption 
$f(l)$ bounds a disk which does not intersect 
$B(f(x_0), kR+kC+3C)$.

Let  $y_{1},...,y_t$ be the vertices of the simplicial complex 
$\varphi(D^2)$ bounded by  the loop $f(l)$. The  
vertices $f(x_1),...,f(x_n)$ are contained among the  $y_j$'s. 
One can suppose that any triangle  $[y_i,y_j,y_m]$ of $\varphi(D^2)$ has   
edge length at most $d$. 
Therefore  we have:  
\[ d(g(y_j),x_i)\leq d(g(y_j),gf(x_i))+C \leq kd(y_j,f(x_i))+2C \leq k^2\varepsilon+(k+2)C.\] 
This proves that $x_i,x_j,g(y_m)$ span a  simplex of $P_a(H)$ 
(for all $i,j,m$) 
whenever we choose $a$ larger than  $k^2 \varepsilon +(k+2)C$. Moreover:  
\[ d(x_0,g(y_i))\geq d(gf(x_0),g(y_i))-C \geq \frac{d(f(x_0),y_i)-3C}{k}-C\geq R.\]  
Further there is a simplicial  map $\psi:\varphi(D^2)\to P_a(H)$ which sends 
$f(x_j)$ into $x_j$ and all other vertices $y_k$ into the
corresponding $g(y_k)$. It is immediate now that 
$\psi\varphi(D^2)$ is a simplicial sub-complex bounded by $l$, 
which has the required properties.   
\end{proof} 
\noindent This proves proposition 1. 
\end{proof}

\noindent When $G$ has torsion, one can construct a highly connected polyhedron with a  free and cocompact $G$-action as follows (see \cite{BM}):

\begin{defi} 
The colored Rips complex $P(d,m,G)$ (for natural $m$) is the sub-complex of the $m$-fold join $G\ast G.....\ast G$ consisting of those simplexes whose vertices are at distance  at most $d$ in $G$. 
\end{defi} 
 
\begin{lem} 
For $m\geq 3$ and  $d$ large enough, $G$ acts freely on the 2-skeleton of $P(d,m,G)$ and  
$\pi_1(P(d,m,G))=0$. 
\end{lem} 
\begin{proof}  
Clearly $G$ acts freely on the vertices of 
$P_d(G)$, hence any non-trivial $g\in G$ fixing a simplex 
has to permute its vertices. 
Adding $m\geq 3$ colors prevents therefore  the action from having  
fixed simplexes of dimension less than $3$. 
Now using the proof of lemma 1 one obtains also the simple connectivity. 
\end{proof} 

\noindent The theorem 1 follows now from the 
proof of proposition 1, suitably adapted to the 2-skeleton of $P(d,m,G)$.  
 
\begin{rem} 
 The same technique shows that the higher connectivity at  
infinity is  also a quasi-isometry invariant of groups.  
\end{rem} 
  
\noindent We have then: 

\begin{cor}
The equivalence class of $V_X(r)$ is 
a quasi-isometry invariant of $X$. 
\end{cor}

\begin{proof}
The result is implied by theorem 1 and lemma 3. 
\end{proof}

\section{Uniform lattices in Lie groups}
\begin{prop}
Uniform lattices in (non-compact) semi-simple Lie groups have linear 
rate of vanishing of $\p1i$. 
\end{prop}
\begin{proof}
We will denote below by $d_X$ the distance function and by $B_X$ the 
respective metric balls for the space $X$.  

Let $K$ be the maximal compact subgroup of the simple Lie group $G$ and
$G/K$ the associated  symmetric space.     
It is well-known that the Killing metric on $G/K$ is  
non-positively curved, and hence the metric balls are diffeomorphic
to standard balls, by the Hadamard theorem. 

If $G$ is not $SL(2,{\bf R})$ then $K$ is 
different from $S^1$ and therefore it has finite fundamental group. 
In particular the universal covering $\widetilde{K}$ is compact. 
The Iwasawa decomposition $G=K A N$ yields a canonical 
diffeomorphism $G\to K\times G/K$. Furthermore we have a  
an induced  canonical quasi-isometry 
$\widetilde{G}\to  \widetilde{K}\times G/K$.
Large balls in $\widetilde{G}$ can be therefore compared with 
products of the (compact) $\widetilde{K}$ and metric balls 
in $G/K$, as follows: 
\[ \widetilde{K}\times B_{G/K}\left({r}-{a}\right)\subset 
B_{\widetilde{G}}(r)\subset \widetilde{K}\times B_{G/K}(a+r)\subset 
B_{\widetilde{G}}(2a+r), 
\mbox{ for } r \mbox{ large enough}, \]
which implies our claim.

The case of $V_{\widetilde{SL(2,{\bf R})}}$ is quite similar, and 
a consequence of the (well-known) fact that 
$\widetilde{SL(2,{\bf R})}$ and $H^2\times {\bf R}$
are canonically  quasi-isometric. We will sketch a proof below. 
Let us outline the construction of the $\widetilde{SL(2,{\bf R})}$ geometry. 
A Riemannian metric on a manifold allows us to construct canonically 
a Riemannian metric on its tangent bundle, usually called the Sasaki 
metric. In particular one considers the restriction of 
the Sasaki metric to the unit tangent bundle $UH^2$ of the 
hyperbolic plane. 
Further there exists a natural diffeomorphism between 
$UH^2$ and  $PSL(2,{\bf R})$, which gives $PSL(2,{\bf R})$ 
a Riemannian metric, and hence induces a metric on its universal covering, 
namely $\widetilde{SL(2,{\bf R})}$. This is the Riemannian 
structure describing the geometry of $\widetilde{SL(2,{\bf R})}$. 
Observe that $\widetilde{SL(2,{\bf R})}$
and $H^2\times {\bf R}$ are two metric structures 
on the same manifold, and both are Riemannian fibrations over $H^2$
(the former being metrically non-trivial while the later is trivial).  

It is clear now that the identity map between the manifolds
$UH^2$ and $H^2\times S^1$ is a quasi-isometry, lifting to a 
quasi-isometry between $\widetilde{SL(2,{\bf R})}$ and $H^2\times {\bf
  R}$.
 This implies that there are two constants $a >0, b$ such that:
\[ \frac{1}{a} d_{H^2\times {\bf R}}(x,y) - b\leq d_{\widetilde{SL(2,{\bf R})}}(x,y) \leq 
a d_{H^2\times {\bf R}}(x,y)+b, \mbox{ for all } x,y\in H^2\times {\bf R}, \]
holds true. 
In particular we have the following inclusions between the 
respective metric balls:
\[ B_{H^2\times {\bf R}}\left(\frac{r}{c}\right)
\subset B_{\widetilde{SL(2,{\bf R})}}(r)\subset B_{H^2\times {\bf R}}(cr)
\subset B_{\widetilde{SL(2,{\bf R})}}(c^2r),  \]
for $r$ large enough and $c >0$. The claim follows. 
\end{proof}
\begin{rem}
The same argument shows that the acyclicity radius for semisimple 
Lie groups is linear (see \cite{Gr2}, section 4).  
\end{rem}

The way to prove the claim for nilpotent and solvable groups consists
in the large scale comparison with some other metrics, 
whose balls are known to be diffeomorphic to standard balls. 
While locally the Riemannian geometry of  a nilpotent Lie group
is Euclidean, globally it is similar to the Carnot-Caratheodory
non-isotropic geometry.  

\begin{prop}
If $G$ is a torsion-free nilpotent group then $V_G$ is linear. 
\end{prop}
\begin{proof}
It is known (see \cite{Ma}) that $G$ is a cocompact lattice in a 
real simply connected nilpotent Lie group $G_{\bf R}$,
called the Malcev completion of $G$.  
We have also a nice  characterization  of the metric
balls in real, nilpotent Lie groups given by Karidi (see \cite{Ka}), 
as follows. Since $G_{\bf R}$ is 
diffeomorphic to ${\bf R}^n$ it makes sense to talk about 
parallelepipeds with respect to the usual Euclidean structure 
on ${\bf R}^n$. Next, there exists some constant
$a > 0$ (depending on the group $G_{\bf R}$ and on the left invariant 
Riemannian structure chosen, but  not on the 
radius $r$) such that the radius $r$-balls $B_{G_{\bf R}}(r)$ are sandwiched 
between two parallelepipeds, which are homothetic at ratio $a$, so:
\[ P_r \subset B_{G_{\bf R}}(r) \subset a P_r\subset B_{G_{\bf R}}(ar),
\mbox{ for any } r \geq 1. \]
This implies that we can take 
$V_{G_{\bf R}}(r) \sim a r$, and hence $V_G$ is linear. 
\end{proof}

\begin{rem} 
Metric balls  in solvable Lie groups are quasi-isometric 
with those of discrete solvgroups, and so they are   
highly concave (see \cite{F}): there
exist pairs of points at distance $c$,   
sitting on the sphere of radius $r$,
which cannot be connected by a path within the ball of radius $r$, 
unless its length is at least $r^{0.9}$. Further it can be shown that 
there are arbitrarily large metric balls which are not simply
connected. Nevertheless we will prove that metric balls 
contain large slices of hyperbolic balls.  
\end{rem}

\begin{prop}
Cocompact lattices in 
solvable stabilizers of horospheres in product of symmetric
  spaces of rank at least two, or generic horospheres in products of
  rank one symmetric spaces have linear $V_G$. 
\end{prop}
\begin{proof}
We give the proof for our favorite solvable group, namely the
3-dimensional group {\sf Sol}.  
It is well-known that {\sf Sol} is isometric to a generic
horosphere ${\mathcal H}$ in the product $H^2\times H^2$ of two 
hyperbolic planes. Generic means here that the horosphere is
associated to a geodesic ray which is neither vertical nor horizontal.
The argument in (\cite{Gr2} 3.D.), or its generalization from  
\cite{Dr}, shows that such horospheres ${\mathcal H}$ are undistorted 
in the ambient space i.e. there exists $a\geq 1$, such that 
\[ \frac{1}{a} d_{H^2\times H^2}(x,y) \leq d_{\mathcal H}(x,y) \leq 
a d_{H^2\times H^2}(x,y), \mbox{ for all } x,y\in {\mathcal H}, \]
holds true. Here $d_{H^2\times H^2}$ and $d_{\mathcal H}$ denote 
the distance functions in $H^2\times H^2$ and ${\mathcal H}$,
respectively. In particular we have the following inclusions between the 
respective metric balls:
\[ {\mathcal H} \cap B_{H^2\times H^2}\left(\frac{r}{a}\right)
\subset B_{\sf Sol}(r)\subset {\mathcal H} \cap B_{H^2\times H^2}(ar)
\subset B_{\sf Sol}(a^2r). \]
Since the horoballs in $H^2\times H^2$ are convex it follows that 
the intersections ${\mathcal H} \cap B_{H^2\times H^2}(r)$ are 
diffeomorphic to standard balls. This proves that 
$V_{\sf Sol}$ is linear. 

The linearity result extends without modifications 
to lattices in solvable stabilizers of generic horospheres 
in symmetric spaces of rank at least 2 (see \cite{Dr}).
\end{proof}
\begin{rem} 
It is known that finitely presented solvable groups are 
either simply connected at infinity or are of a very special form, 
as described in \cite{Mi}. On the other hand it is a classical 
result that any simply connected solvable Lie group is diffeomorphic 
to the Euclidean space. It would be interesting to know whether 
a simply connected solvable Lie group can be isometrically 
embedded as a horosphere in a symmetric space.      
\end{rem}
\begin{rem}
Notice that horospheres in hyperbolic spaces (and hence non-generic 
horospheres in products of hyperbolic spaces) have exponential
distortion, namely  $d_{H^n}(x,y)\sim \log \, d_{{\mathcal H}^{n-1}}(x,y)$, for 
$x,y \in {\mathcal H}^{n-1}$. This highly contrast with the higher rank
and/or generic case. 
\end{rem}
\noindent This ends the proof of theorem 2.

\begin{rem}
One might notice a few similarities between $V_G$ and the isodiametric
function considered by Gersten (see \cite{Gr2}).
\end{rem}

\begin{cor}
The rate of  vanishing of $\p1i$ is linear for the fundamental 
groups of geometric 3-manifolds. 
\end{cor}
\begin{proof}
There are eight geometries in the 
Thurston classification (see  \cite{Sc}): the sphere $S^3$, 
$S^2\times {\bf R}$, the Euclidean $E^3$, the hyperbolic 3-space $H^3$, 
$H^2\times {\bf R}$, $\widetilde{SL(2,{\bf R})}$, {\sf Nil} and 
{\sf Sol}. Manifolds covered by  $S^3$ have finite fundamental 
groups. Further the compact manifolds without boundary covered 
by $S^2\times {\bf R}$ are the two $S^2$ bundles over $S^1$, 
${\bf R P}^2\times S^1$ or  the connected sum 
${\bf R P}^3\sharp {\bf R P}^3$, and the claim can be checked 
easily. As we already observed, this is the also the case  for the 
Euclidean and hyperbolic geometries.  The same holds 
for the product $H^2\times {\bf R}$, in which case metric balls 
are diffeomorphic to standard balls. The remaining cases are covered
by theorem 2.  
\end{proof}

\bibliographystyle{plain}

\begin{thebibliography}{10} 
 
\bibitem{BM} 
M.Bestvina and G.Mess, {\normalfont\itshape The boundary of negatively
  curved groups}, J.A.M.S. 
{\bf{4}}(1991), 469-481. 
 
\bibitem{Br} 
S.Brick, {\normalfont\itshape Quasi-isometries and ends of groups}, J.Pure Appl. Algebra {\bf{86}}(1993), 23-33. 
 
\bibitem{BrM} 
S.Brick and M.Mihalik, {\em The QSF property for groups and spaces}, Math. Zeitschrift  
{\bf 220}(1995), 207-217.  
  

\bibitem{D}
M.Davis, {\em Groups generated by reflections and aspherical manifolds
  not covered by Euclidean space}, Ann. of Math.  {\bf 117}(1983), 293--324.

\bibitem{DE}
A.Dyubina-Erschler, {\em Instability of the virtual solvability and the
  property of being virtually torsion-free for quasi-isometric
  groups},  I.M.R.N. 2000, no. 21, 1097--1101.

\bibitem{Dr}
C.Dru\c tu, {\em Nondistorsion des horosph\`eres dans des immeubles 
euclidiens et dans des espaces sym\'etriques}, G.A.F.A. {\bf 7}(1997), 712--754.


\bibitem{F}
L.Funar, {\em Discrete cocompact solvgroups and Po\'enaru's condition}, 
Archiv Math.(Basel) {\bf 72}(1999), 81-85.

\bibitem{GdlH} 
E.Ghys and P.de la Harpe (Ed.), {\normalfont\itshape Sur les groupes
  hyperboliques d'apr\`es M. Gromov}, Progress in Math.,
vol. {\bf{3}}, Birkhauser, 1990. 
 
\bibitem{GM} 
R.Geoghegan and M.Mihalik, {\normalfont\itshape The fundamental group
  at infinity}, Topology  {\bf{35}}(1996), 
  655-669. 
 
%\bibitem{Ge}
%S.Gersten, {\em Isoperimetric and isodiametric functions of finite
%  presentations}, Geometric group theory, Vol. 1 (Sussex, 1991), 79--96, 
%London Math. Soc. Lecture Note Ser., 182, (G.A.Niblo and M.A.Roller, Ed.), 
%Cambridge Univ. Press, Cambridge, 1993.


\bibitem{Gr} 
M.Gromov, {\normalfont\itshape Hyperbolic groups}, Essays in Group Theory (S. Gersten Ed.), MSRI publications, no. {\bf{8}}, Springer-Verlag (1987). 


\bibitem{Gr2}
M.Gromov, {\em  Asymptotic invariants of infinite groups}, 
Geometric group theory, Vol. 2 (Sussex, 1991), 1--295, 
London Math. Soc. Lecture Note Ser., 182, (G.A.Niblo and M.A.Roller, Ed.),
Cambridge Univ. Press, Cambridge, 1993.

\bibitem{Ka}
R.Karidi, {\em Geometry of balls in nilpotent Lie groups}, 
Duke Math. J. {\bf 74}(1994),  301--317. 

\bibitem{Ma}
A.Malcev, {\em On a class of homogeneous spaces}, 
Izvestiya Akad. Nauk. SSSR. Ser. Mat. {\bf 13}(1949), 9--32,  
Amer. Math. Soc. Translation  {\bf 39}(1951).  

\bibitem{Mi}
M.Mihalik, {\em Solvable groups that are simply connected at
  $\infty$},  Math. Zeitschrift {\bf 195} (1987), 79--87. 


\bibitem{Ot} 
D.E.Otera, {\normalfont\itshape On the simple connectivity of groups},
preprint 139(2001), Univ. Palermo, Bull.U.M.I. (to appear). 
 

\bibitem{Sc}
P.Scott, {\em The geometries of 3-manifolds}, Bull.London Math. Soc.,
{\bf 15}(1983), 401-487.  

\bibitem{Ta} 
 
C.Tanasi, {\normalfont\itshape Groups simply connected at infinity}, 
(Italian), Rend. Ist. Mat. Univ. Trieste   
{\bf{31}} (1999), 61-78. 
 


\end{thebibliography}

\end{document}